\documentclass[12pt,thmsa]{article}
\usepackage{amsmath}
\usepackage{amssymb}

\begin{document}

\centerline{\bf IMAGE MEASURES OF INFINITE PRODUCT MEASURES  }

\centerline{\bf AND  GENERALIZED BERNOULLI CONVOLUTIONS }

\vskip 10mm

\centerline {SERGIO ALBEVERIO}

\vskip 3mm

\centerline {Institut f\"{u}r Angewandte Mathematik,
Universit\"{a}t Bonn, }

\centerline{  Wegelerstr. 6, D-53115 Bonn (Germany);}

\centerline{ SFB 256, Bonn,  BiBoS, (Bielefeld - Bonn); }

\centerline{ IZKS Bonn; CERFIM, Locarno and Acc. Arch. (USI)
(Switzerland).}

\vskip 5mm

\centerline {GRYGORIY TORBIN}

\vskip 3mm

\centerline {National Pedagogical University, Kyiv (Ukraine)}

\vskip 5mm

\begin{abstract}

\thinspace We examine measure preserving mappings $f$ acting from a
probability space $\left( \Omega ,F,\mu \right) $ into a probability space $%
\left( \Omega ^{*},F^{*},\mu ^{*}\right) ,$ where $\mu ^{*}=\mu
(f^{-1})$. Conditions on $f$, under which $f$ preserves the
relations ''to be singular''\thinspace and ''to be absolutely
continuous''\thinspace between measures defined on $\left( \Omega
,F\right) $ and corresponding image measures, are investigated.

We apply the results to investigate the distribution of the random variable $%
\xi =\sum\limits^{\infty}_{k=1} \xi _k\lambda ^k,$ where $%
\lambda \in (0;1),$ and $\xi _k$ are independent not necessarily identically
distributed  random variables taking the values $i$ with probabilities $%
p_{ik}$ ,$i=0,1.$

We also studied in details the metric-topological  and  fractal
properties of the distribution of a random variable $\psi =
\sum\limits^{\infty }_{k=1} \xi _ka_k,$ where $a_k>0$ are terms of
the convergent series.

\end{abstract}

\smallskip \smallskip $^1$I{nstitut f\"{u}r Angewandte Mathematik,
Universit\"{a}t Bonn, Wegelerstr. 6, D-53115 Bonn\ (Germany) \thinspace
\thinspace \thinspace }$^2${SFB 611, \ Bonn, \ BiBoS, Bielefeld - Bonn
\thinspace \thinspace \thinspace }$^3${CERFIM, Locarno and USI
(Switzerland)\thinspace \thinspace \thinspace \thinspace }$^4${\ IZKS, Bonn}

$^5${National Pedagogical University}, Kyiv (Ukraine) {E-mail:
torbin@imath.kiev.ua}

\smallskip \smallskip

\textbf{AMS Subject Classifications (2000): 11K55,\thinspace \thinspace
28A35, 28A80,\thinspace 30B20, 60E10, \thinspace 60G30.}\medskip

\textbf{Key words:} Bernoulli convolutions, Hausdorff-Besicovitch dimension,
fractals, singular probability distributions, absolutely continuous
probability distributions, nonlinear projections, image measures.

\section{Introduction}

Let $P_{\,\xi }$ be the probability distribution of a random variable
\begin{equation}
\xi =\sum^{\infty }_{k=0} \xi _k\lambda ^k, \label{Ksi}
\end{equation}
where $\lambda \in (0;1)$, and $\xi _k$ are independent random
variables taking values $\,-1$ and $+1$ with probabilities $\frac
12$ . $\,$ $P_{\,\xi }$ is known as ''infinite Bernoulli
convolution''. Measures of this form have been studied since
1930's from the pure probabilistic point of view as well as for
their applications in harmonic analysis, dynamical systems and
fractal analysis. We will not describe in details the history of
the investigation of these measures, paper \cite{Peres} contains a
comprehensive survey on Bernoulli convolutions, corresponding
historical notes and brief discussion of some applications,
generalizations and problems.

Let $S_{\perp }$ be the set of those $\lambda \in [\frac 12;1)$ for which
the random variable $\xi $ has singular distribution, and let $S_{\ll }$ be
the set of those $\lambda \in (\frac 12;1)$ for which the measure $P_{\,\xi
} $ is absolutely continuous (with respect to Lebesgue measure $L$%
).\thinspace In 1939 P. Erd\"{o}s \cite{Erdos} proved that $\lambda \in $ $%
S_{\perp }$whenever $\lambda $ is the reciprocal of a Pisot number in $(1;2).
$ Up to now we do not know any other examples of $\lambda $ belonging to $%
S_{\perp }. $ In 1940 P. Erd\"{o}s proved the existence a number $a$
independent of $\lambda $ such that almost all $\lambda $ from $(a;1)$
belong to $S_{\ll }.$ In 1962 A. Garsia \cite{Garsia}\thinspace found the
largest explicitly described subset of $S_{\ll }$ known up to now. Garsia
formulated the following conjecture: almost all $\lambda \in [\frac 12;1)$
belong to $S_{\ll }.$ This conjecture was proved by B.Solomyak\cite{Sol1} in
1995, a simpler proof was found by Y. Peres and B.Solomyak in \cite{PeresSol}%
. In \cite{PeresSolom2} a more general nonsymmetric model was considered by
Y. Peres and B.Solomyak: if the $\xi _k$ in (\ref{Ksi}) take values $-1$ and $%
+1$ with probabilities $1-p$ and $p$ ($p\in [\frac 13;\frac 23])$%
correspondingly, then almost all $\lambda $ from $[p^p(1-p)^{(1-p)};1)$
belong to $S_{\ll }.$

\smallskip The main goal of this paper is to investigate the
distributions of the following random variables which are generalizations of
the Benoulli convolutions:

1)$\,\psi = \sum\limits^{\infty }_{k=1}\psi _k\lambda ^k,$%
where $\lambda \in (0;1)$, and $\psi _k$ are independent not necessarily
indentically distributed random variables taking the values $\,0$ resp.$1$
with probabilities $p_{0k}$ resp. $p_{1k}$ , i.e., schematicall\.{y}:
\[
\begin{array}{lll}
\psi _k & 0 & 1 \\
& p_{0k} & p_{1k}
\end{array}
,
\]
with $p_{ik}\geq 0,\,p_{0k}+p_{1k}=1.$

2)$\varphi =\sum\limits^{\infty }_{k=1} \varphi _ka_k,$where $%
\varphi _k$ are independent random variables with the following
distributions
\[
\begin{array}{lll}
\varphi _k & 0 & 1 \\
& p_{0k} & p_{1k}
\end{array}
,
\]
$p_{ik}\geq 0,\,p_{0k}+p_{1k}=1$, and $a_k\geq 0$ are terms of a
convergence series $\sum\limits^{\infty }_{k=1} a_k.$

The main method we use is the method of image measures which are nonlinear
projections of infinite product measures. In section 2 we consider an
infinite product measure spaces and the problem of absolute continuity,
singularity and discreteness of the product measures. Section 3 is devoted
to measure preserving mappings $f$ acting from a probability space $\left(
\Omega ,F,\mu \right) $into a probability space $\left( \Omega
^{*},F^{*},\mu ^{*}\right) ,$ where $\mu ^{*}=\mu (f^{-1})$. Conditions on $%
f $, under which $f$ preserves the relations ''to be singular''\thinspace
and ''to be absolutely continuous''\thinspace between measures defined on $%
\left( \Omega ,F\right) $ and corresponding image measures, are investigated.

In section 4, by using the results of sections 2 and 3, we prove a slightly
sharpened versions of the results of B.Solomyak\cite{Sol1}\thinspace ,and
Y.Peres and B. Solomyak \cite{PeresSolom2}.\thinspace We also consider an
counterexample concerning a problem which was formulated in section 3.

Section 5 is devoted to the study of metric, fractal and topological
properties of the distribution of the above mentioned random variable $%
\varphi .$ We prove criteria for singularity and absolutely
continuity of the distribution of the r.v. $\varphi $ in the case,
when $a_k\geq \sum\limits^\infty_{i=k+1} a_i$.

\section{Infinite product measure spaces and absolutely continuity of
product measures.}

Let us consider an independent sequence of probability spaces $\left\{
\left( \Omega _k,A_k,\mu _k\right) \right\} $. We will denote by $\left(
\Omega ,A,\mu \right) $ the infinite product of the probability spaces $%
\left( \Omega _k,A_k,\mu _k\right) :$%
\[
\left( \Omega ,A,\mu \right) =\prod^{\infty }_{k=1} \left( \Omega
_k,A_k,\mu _k\right) .
\]
\thinspace For any element $\omega =\left( \omega _1,\omega
_2,...,\omega _k,...\right) \in \Omega \,$ we have
\[
\mu (\omega )=\prod^{\infty }_{k=1} \mu _k(\omega
_k),\,\,\,\,\omega _k\in \Omega _k.
\]

\textbf{Theorem 1.}$\,\mu \,$\textit{is pure discrete if and only if }
\begin{equation}
\prod^{\infty }_{k=1} \max_{\omega _k\in \Omega _k} \mu _k(\omega
_k)>0. \label{DiscrMu}
\end{equation}

\smallskip
\textbf{Proof.} If $\prod\limits^{\infty }_{k=1}
\max\limits_{\omega _k\in \Omega _k} \,\mu _k(\omega _k)=0,$ then
for any point $\omega \in \Omega $ we have $\mu (\omega )\leq $
$\prod\limits^{\infty }_{k=1} \max\limits_{i\in \Omega _k}
p_{ik}=0.$ Therefore, condition \ref{DiscrMu} is necessary for
discreteness of the measure $\mu.$

To prove the sufficiency we consider a subset $A_{+}\subset \Omega :$%
\[
A_{+}=\left\{ \omega :\,\mu _k(\omega _k)>0\text{ and }\,
\prod^{\infty }_{k=1} \mu _k(\omega _k)>0\right\} .
\]

The set $A_{+}$ consists of the points $\omega =\left( \omega
_1,\omega _2,...,\omega _k,...\right) $ such that $\,\mu _k(\omega
_k)>0$ and condition $\mu _k(\omega _k)\neq \, \max\limits_{\omega
_k\in \Omega _k} \,\mu _k(\omega _k)$ holds only for finite
numbers of values $k.$ It easy to see that the set $A_{+}\,$is at
most countable and the event ''$\omega \in A_{+}$'' does not
depend on any finite coordinates of $\omega .$ Therefore, by using
the ''0 and 1'' theorem of Kolmogorov, we conclude that $\mu
(A_{+})=0\,$ or $\mu (A_{+})=1.$ Since the set $A_{+}$ contains the point $%
\omega ^{*}\,$ such that $\mu _k(\omega
_k^{*})=\max\limits_{\omega _k\in \Omega _k} \,\mu _k(\omega _k)$,
we have $\mu (A_{+})\geq \mu (\omega ^{*})>0.$ Thus, $\,\mu
(A_{+})=1,$ which proves the discreteness of the measure $\mu .$
$\square$

Let $\left( \Omega ,A,\nu \right) = \prod\limits^{\infty }_{k=1}
\left( \Omega _k,A_k,\nu _k\right) $ such that the measure $\nu
_k\;$ is absolutely continuous with respect to the measure $\mu
_k\,$ (for short $\nu _k\ll \mu _k$) for all $k\in N.$ By using
completely analog arguments as in \cite{Kak},\thinspace one can
prove the following sharper variant of the Kakutani theorem.

\textbf{Theorem 2 }\textit{Let }$\nu _k\;$\textit{be absolutely continuous
with respect to the measure }$\mu _k\,$\textit{. Then the measure }$\mu $%
\textit{\ is either purely absolutely continuous with respect to the measure
}$\mu $\textit{\ or purely singular (including the discreteness). Moreover, }
\begin{equation}
\nu \ll \mu \,\,\text{if and only if } \prod^{\infty }_{k=1} \rho
(\mu _k,\nu _k)>0,  \label{AbsContKAK}
\end{equation}

\textit{where }$\rho (\mu _k,\nu _k)= \int\limits_{\Omega _k} \sqrt{%
\frac{d\nu _k}{d\mu _k}}d\mu _k.$

\textbf{Remark. }
 The expression $\int\limits_{\Omega _k} \sqrt{\frac{d\nu _k%
}{d\mu _k}}d\mu _k\,$coincides with the Hellingers integral \cite{Kak}.

\section{Measure-preserving mappings of probability spaces.}

Let \smallskip $\left( \Omega ,F,\eta \right) $and $\left( \Omega ,F,\tau
\right) $ be abstract probability spaces.

Let us consider a measurable mapping $f:$
\begin{eqnarray*}
&&(\Omega ,F,\eta )\stackrel{f}{\rightarrow }(\Omega ^{*},F^{*},\eta ^{*}),
\\
&&\left( \Omega ,F,\tau \right) \stackrel{f}{\rightarrow }(\Omega
^{*},F^{*},\tau ^{*}),
\end{eqnarray*}
\smallskip

where the measures $\eta ^{*}$ and $\tau ^{*}\,$are defined as follows:
\begin{eqnarray*}
\,\,\eta ^{*}(E) &=&\,\eta \,(f^{\,\,-1}(E))\,, \\
\,\tau ^{*}(E) &=&\,\tau \,(f^{\,\,-1}(E))\,,\,\,\text{for any subset
\thinspace }E\in F^{*}\,,\,\,\,\,\,
\end{eqnarray*}

\smallskip and $f^{\,\,-1}(E)=\left\{ \omega :\,\omega \in \Omega \,\text{and%
}\,\,\,f(\omega )\in E\right\} .$

The mapping $f\,\,$is measure-preserving by definition (for details see
,e.g., \cite{Bil},\cite{WAL}).

\textbf{Theorem 3.} \textit{If }$\eta \ll \tau ,$\textit{\ then }$\eta
^{*}\ll \tau ^{*}.$

\textbf{Proof.} \smallskip Suppose that $\,\,\tau ^{*}(E)=0$ for
some subset $E\in F^{*}.$ Then $\tau (f^{-1}(E))=\tau ^{*}(E)=0.$
Since $\eta \ll \tau ,$we conclude that $\eta (f^{-1}(E))=\eta
^{*}(E)=0,$ which proves an absolute continuity
of the measure $\eta ^{*}\,$ with respect to the measure$\,\,\tau ^{*}$.%
$\square$

\textbf{Remark }In general the inverese result to theorem 3 does
not hold. In fact, it is possible to construct probability spaces
\smallskip $\left( \Omega ,F,\eta \right) $and $\left( \Omega
,F,\tau \right) ,$and a mapping $f $ such that $\eta \perp \tau
\,,$but $\eta ^{*}\ll \tau ^{*}.$ A corresponding example will be
considered later in section 4.

\textbf{Theorem 4.} \textit{If }$\eta ^{*}\perp \tau ^{*},$\textit{then }$%
\eta \perp \tau .$

\textbf{Proof.} If $\eta ^{*}\perp \tau ^{*},$then there exists a subset $E$ $\in F^{*}\,$%
such that $\eta ^{*}(E)=0$ and $\tau ^{*}(\Omega ^{*}\backslash
E)=0.\,$By the definitions of the measures $\eta ^{*}$ and $\,\tau
^{*},$  we have $\eta \,(f^{\,\,-1}(E))=\eta ^{*}(E)=0$ and
$\,\tau (f^{\,\,-1}(\Omega
^{*}\backslash E))=\tau ^{*}(\Omega ^{*}\backslash E)=0.$ Since $%
f^{\,\,-1}(\Omega ^{*}\backslash E)\cap f^{\,\,-1}(E)=\emptyset ,$
we conclude for the singularity of the measures $\eta \,$and$\,\tau $%
$\square$

\textbf{Theorem 5. }\textit{If }$f$\textit{\thinspace is a bijective
mapping, then }
\begin{eqnarray*}
\eta &\ll &\tau \,\,\text{if and only if }\eta ^{*}\ll \tau ^{*}\,\text{ },
\\
\eta &\perp &\tau \,\,\text{if and only if }\eta ^{*}\perp \tau ^{*}.
\end{eqnarray*}

\textbf{Proof.} If $f$ \thinspace is a bijective mapping, then
$A=f^{\,-1}[f(A)]$ and $\eta (A)=\eta ^{*}(f(A)).$ Using the same
arguments, we have $\tau (A)=\tau ^{*}(f(A))$ for all measurable
subsets $A\subset \Omega .$

To prove the first statement of the theorem it is sufficient to
prove that the condition $\eta ^{*}\ll \tau ^{*}\,$implies the
condition $\eta \ll \tau .\,$ Suppose that $\,\,\tau (A)=0$ for
some measurable subset $A\subset
\Omega .$ Then $\tau ^{*}(f(A))=\tau (A)=0.$ Since $\eta ^{*}\ll \tau ^{*},$%
we conclude that $\eta ^{*}(f(A))=\eta (A)=0,$ which proves the absolute
continuity of the measure $\eta \,$ with respect to the measure$\,\,\tau .$

\smallskip The second assertion of the theorem can be prove in a similar way.%
$\square$

\textbf{Theorem 6.\thinspace }\textit{If there exists a measurable subset }$%
\Omega _0\subset \Omega \,$\textit{such that \smallskip }$\eta (\Omega
_0)=\tau (\Omega _0)=0$\textit{\ and the mapping }$f:$\textit{\ }$(\Omega
\backslash \Omega _0)\rightarrow \Omega ^{*}$\textit{\ }$\,$\textit{is a
bijection, then }
\begin{eqnarray*}
\eta &\ll &\tau \,\,\text{if and only if }\eta ^{*}\ll \tau ^{*}\,\text{ },
\\
\eta &\perp &\tau \,\,\text{if and only if }\eta ^{*}\perp \tau ^{*}.
\end{eqnarray*}

\smallskip

\textbf{Proof.} For any measurable subset $A\subset \Omega \,$ we
have $A\subset f^{\,-1}[f(A)].$ Since $f:(\Omega \backslash \Omega
_0)\rightarrow \left[ 0;1\right] $ is a bijective mapping, we
conclude that $f^{\,-1}[f(A)]=A\cup A_0,$ where $A_0\subset \Omega
_0.$ Therefore,
\begin{eqnarray*}
\eta (A) &=&\eta (A\cup A_0)=\eta (f^{\,-1}[f(A)])=\eta ^{*}(f(A)), \\
\text{and \thinspace \thinspace \thinspace \thinspace \thinspace }\tau (A)
&=&\tau ^{*}(f(A)).
\end{eqnarray*}

We complete the proof by arguments which are completely similar to
those used in the proof of the previous theorem. $\square$

\section{Generalized Bernoulli convolutions I}

Let us consider a random variables
\[
\psi =\sum^{\infty }_{k=1} \psi _k\lambda ^k,
\]
where $\lambda \in (0;1)$, and $\psi _k$ are independent not necessarily
indentically distributed random variables with the following distributions
\[
\begin{array}{lll}
\psi _k & 0 & 1 \\
& p_{0k} & p_{1k}
\end{array}
,
\]
with $p_{ik}\geq 0,\,p_{0k}+p_{1k}=1.$

\smallskip

\textbf{Theorem 7.} \textit{If }
\begin{equation}
\sum^{\infty }_{k=1} \left( \frac 12-p_{0k}\right) ^2<\infty ,
\label{Sum1/2}
\end{equation}
\textit{then for almost all }$\lambda \in [\frac 12;1)$\textit{\ the r.v. }$%
\psi $\textit{\ has an absolute continuous distribution.}

\textbf{Proof.} Let us consider two sequences of probability
spaces: $\left\{ (\Omega _k,A_k,\mu _k)\right\} $ and $\left\{
(\Omega _k,A_k,\nu _k)\right\} ,$
where $\Omega _k=\left\{ 0;1\right\} ,$ $A_k$ consists of all subsets of $%
\Omega _k$ and measures $\mu _k$ and $\nu _k$ are defined as follows:
\begin{eqnarray*}
\nu _k(0) &=&p_{0k},\,\,\,\nu _k(1)=p_{1k}; \\
\mu _k(0) &=&\frac 12,\,\,\,\,\,\,\mu _k(1)=\frac 12.
\end{eqnarray*}
It is easy to see that $\nu _k\ll \mu _k$ for all $k$, and $\rho (\mu _k,\nu
_k)=\int\limits_{\Omega _k} \sqrt{\frac{d\nu _k}{d\mu _k}}d\mu _k=\sqrt{%
\frac 12p_{0k}}+\sqrt{\frac 12p_{1k}}.$ Let $\left( \Omega ,A,\mu \right) =%
\prod\limits^{\infty }_{k=1} \left( \Omega _k,A_k,\mu _k\right)
$and $\left( \Omega ,A,\nu \right) = \prod\limits^{\infty }_{ k=1}
\left( \Omega _k,A_k,\nu _k\right) .$ By using theorem 2, we can
deduce that the measure $\nu $ is absolute continuous with respect
to the measure $\mu \,$ if and only if
\begin{equation}
\prod^{\infty }_{k=1} \left( \sqrt{\frac 12p_{0k}}+%
\sqrt{\frac 12p_{1k}}\right) <\infty .  \label{InfProdBern}
\end{equation}
It is easy to check that the product (\ref{InfProdBern}) converges if and
only if condition (\ref{Sum1/2}) holds.

Let us consider a mapping $f$ : $\Omega \rightarrow [0;\frac \lambda
{1-\lambda }]$ defined as follows: for any $\omega =(\omega _1\omega
_2\ldots \omega _k\ldots )\in \Omega ,$ $\,$%
\begin{equation}
f(\omega )= \sum^{\infty }_{k=1} \omega _k\lambda ^k. \label{f}
\end{equation}
It is easy to see that $f$ is a measurable mapping. The image
measure $\mu ^{*}=\mu (f^{-1})\,$is the above classical mentioned
Bernoulli measure $\mu _{\,\lambda },$ the image measure $\nu
^{*}=\nu (f^{-1})$ is the probability measure corresponding to the
distribution of the random variable $\psi .$ By using theorem 2
and theorem 3, we conclude that condition (\ref{Sum1/2}) implies
the absolute continuity of the measure $\nu ^{*}$ with respect to
measure $\mu ^{*}.\,$B.Solomyak in \cite{Sol1} proved that for almost all $%
\lambda \in [\frac 12;1)$ the measure $\mu ^{*}$ is absolutely continuous
with respect to Lebesgue measure. Moreover, since for $\lambda \in [\frac
12;1)$ the support of the measure $\mu ^{*}$ coincides with the whole closed
interval $[0;\frac \lambda {1-\lambda }],$we conclude that $\mu ^{*}$ is
equivalent to Lebesgue measure ( for short $\mu ^{*}\sim L$) for almost all $%
\lambda \in [\frac 12;1)$. Therefore, $\nu ^{*}\ll L$ for almost all $%
\lambda \in [\frac 12;1).$%
$\square$

\textbf{Theorem 8. }\smallskip \textit{If there exist a number
}$p\in [\frac 13;\frac 23]$\textit{\ such that }
\[
\sum^{\infty }_{k=1} \left( p_0-p_{0k}\right) ^2<\infty ,
\]
\textit{then for almost all }$\lambda \in [p^p\cdot (1-p)^{(1-p)}\,;1)$%
\textit{\ the r.v. }$\psi $\textit{\ has an absolutely continuous
distribution.}

\textbf{Proof.} The proof of this theorem is analogous to the
previous one, but we define the measures $\mu _k$ as follows:
\[
\mu _k(0)=1-p,\,\,\,\,\,\,\mu _k(1)=p,
\]
and we use the following result from \cite{Peres}: if $p\in [\frac 13;\frac
23],$ then for almost all $\lambda \in [p^p\cdot (1-p)^{(1-p)}\,;1)$ the$\,$%
measure $\mu $ is absolutely continuous with respect to Lebesgue measure.%
$\square$

\textbf{Proposition} \textit{There exist product measures }$\nu $\textit{\
and }$\mu $\textit{, and a measure preserving mapping }$f$\textit{\ such
that }$\nu $\textit{\ }$\perp $\textit{\ }$\mu ,$\textit{\ but }$\nu ^{*}\ll
$\textit{\ }$\mu ^{*}.$

\textbf{Proof.} Let $\Omega _k=\left\{ 0;1\right\} ,$ $A_k$ consists of all subsets of $%
\Omega _k$ and measures $\mu _k$ and $\nu _k$ defined as follows:
\begin{eqnarray*}
\nu _k(0) &=&p_{0k}=1-p,\,\,\,\,\,\,\,\,\,\nu _k(1)=p_{1k}=p\in
\left[ \frac 13;\frac 23 \right],\,p\neq \frac 12; \\ \mu _k(0)
&=&\frac 12,\,\,\,\,\,\,\,\,\,\,\,\mu _k(1)=\frac 12.
\end{eqnarray*}
Since $\sum\limits^{\infty }_{k=1} \left( \frac
12-p_{0k}\right) ^2=\infty ,$ the corresponding product measures $\nu $ and $%
\mu $ are mutually singular. Let the mapping $f\;$be defined by (\ref{f}).
As mentioned above, the image measure $\mu ^{*}=\mu (f^{-1})\,$is equivalent
to Lebesgue measure for almost all $\lambda \in [\frac 12;1)$, and the image
measure $\nu ^{*}=\nu (f^{-1})$ is equivalent to Lebesgue measure for almost
all $\lambda \in [p^p\cdot (1-p)^{(1-p)}\,;1).$ Therefore, for almost all $%
\lambda \in [p^p\cdot (1-p)^{(1-p)}\,;1)$, the image measures $\nu ^{*}$ and
$\mu ^{*}$ are equivalent.%
$\square$

\section{\protect\smallskip Generalized Bernoulli convolutions II}

Let us consider a random variable
\[
\varphi = \sum^{\infty }_{k=1} \varphi _ka_k,
\]
where $\varphi _k$ are independent random variables with the following
distributions
\[
\begin{array}{lll}
\varphi _k & 0 & 1 \\
& p_{0k} & p_{1k}
\end{array}
,
\]
$p_{ik}\geq 0,\,p_{0k}+p_{1k}=1$, and $a_k\geq 0$ are terms of a convergent
series. Without loss of generality we shall assume that
\begin{equation}
\sum^{\infty }_{k=1} a_k=1. \label{Series1}
\end{equation}
One can proved that the support $S_{\,\varphi }$ of the r.v. $\varphi $,
i.e., the smallest closed set supported the distribution, is a perfect set
of the following form:
\begin{equation}
S_{\,\varphi }=\left\{ x:\,x=\sum^{\infty }_{k=1} \gamma
_k(x)\,a_k,\,\gamma _k(x)\in \left\{ 0;1\right\} \text{
and\thinspace }p_{\gamma _k(x),k}\neq 0\right\} .
\label{Support1}
\end{equation}
In other words, $S_{\,\varphi }$ is a set of all possible
''incomplete sums'' of the series (\ref{Series1}).\thinspace The
metric and topological properties of the set $S_{\,\varphi }$
directly depend on the properties of the series (\ref{Series1}).
Let us consider three examples illustrating this dependence.

\textbf{Example 1.} If $a_k=\frac 2{3^k},$ then $S_{\,\varphi }$ coincides
with the classical Cantor set $C_0.$

\textbf{Example 2}. If $a_k=\frac 1{2^k},$ then $S_{\,\varphi }$ $=[0;1].$

\textbf{Example 3.} If $a_k=\frac{1-\,\varepsilon _0}{2^k}+\frac{%
3\,\varepsilon _0}{4^k},\varepsilon _0\in (0;1),\,$  then, by
using theorem 9, one proves that $S_{\,\varphi }$ is a nowhere
dense set of positive Lebesgue measure.

Let $r_k=$\smallskip $\sum\limits^{\infty }_{i=k+1} a_i$ and
$\delta _k=\frac{a_k}{r_k}.$ If the condition $\delta _k<1$ holds
for all $k>k_{0\text{ }},$then the support $S_{\,\varphi }$ is a
finite union of closed intervals. Moreover, the distribution of
the r.v. $\varphi $ is the so-called distribution with big
overlaps, because almost all points of the support have an
essentially nonunique representation of the form (\ref
{Support1}). The following theorem is proven in (\cite{PrTor}).

\textbf{Theorem 9.} \textit{If condition }$\delta _k>1$\textit{holds for all
}$k>k_0,$\textit{\ then }$S_{\,\varphi }$\textit{\ is a nowhere dense set,
and }$S_{\,\varphi }$\textit{\ has a positive Lebesgue measure if and only
if }
\begin{equation}
\sum^{\infty }_{k=1}(\delta _k-1)<\infty . \label{SerDelta}
\end{equation}
\textit{If condition (\ref{SerDelta}) does not hold, then the r.v. }$\varphi
$\textit{\ has a fractal distribution with the following
Hausdorff-Besicovitch dimension of the support: }
\[
\alpha _0(S_{\,\varphi })=\lim\limits_{\overline{k\rightarrow
\infty }} \frac{k\,\ln 2}{\sum\limits^{k}_{i=1} (\delta _i+1)}.
\]

\textit{In addition, if condition }$\delta _k>1$\textit{\ holds for all }$%
k\in N$\textit{, then all points of the support }$S_{\,\varphi }$\textit{\
have a unique representation of the form (\ref{Support1}).}

Our main goal in this section is to investigate the topological properties
of the support $S_{\,\varphi }$ in a general situation, i.e., when there
exist sequences $\left\{ k_n^{^{\prime }}\right\} $ and $\left\{
k_n^{^{^{\prime }\,^{\prime }}}\right\} $ such that $\delta _{k_n^{^{\prime
}}}>1$ and $\delta _{k_n^{^{\prime \prime }}}<1$ for all $n\in N.$

\textbf{Theorem 10.}\textit{\ Let }$p_{ik}>0$ and \textit{the sequence }$%
\left\{ a_k\right\} \,$\textit{be monotonically nonincreasing, i.e., }$%
a_k\geq a_{k+1},k\in N.$\textit{\ Then the support }$S_{\,\varphi }$\textit{%
\ is a nowhere dense set if and only if condition }$\delta _k>1$\textit{\
holds for infinitely many number of }$k.$

\textbf{Proof.} \textbf{Necessity }is obvious.

\textbf{Sufficiency. }To obtain a contradiction we suppose that the
condition $\delta _k>1$ holds for infinitely many number of $k$, but there
exist an closed interval $[\alpha ^{*};\beta ^{*}]$ such that $[\alpha
^{*};\beta ^{*}]\subset S_{\,\varphi }.$ Let
\[
\alpha =\inf \left\{ x:\,[x;\beta ^{*}]\subset S_{\,\varphi }\right\} ,
\]
and
\[
\beta =\sup \left\{ x:[\alpha ^{*};x]\,\subset S_{\,\varphi }\right\} .
\]
Since $S_{\,\varphi }$ is a closed set, $\alpha \in $ $S_{\,\varphi }$ and $%
\beta \in $ $S_{\,\varphi }.$ Let us prove that $\alpha \neq 0.$ If we
suppose that $\alpha =0$ , then we should conclude that $[0;\beta ]\subset
S_{\,\varphi }.$ We can choose a number $k_1$ such that $a_{k_1}<\beta $ and
$a_{k_1}>r_{k_1.}$ Then, the interval $(r_{k_1};a_{k_1})$ does not contain
any point from the support. Indeed, if $\varphi _1=1$ or $\varphi _2=1$ or$%
\ldots $or $\varphi _{k_1}=1,$ then $\varphi \geq a_{k_1};$ if $\varphi _1=0$
and $\varphi _2=0$ and$\,\ldots $ and $\varphi _{k_1}=0,$then $\varphi \leq
r_{k_1.}$ In a similar way we prove that $\beta \neq 1.$ Since $S_{\,\varphi
}$ is a perfect set, from the construction of the set $[\alpha ;\beta ]\,$%
there follows the existence of a positive number $\varepsilon _0$ such that
both of the intervals $(\alpha -\varepsilon _0;\alpha )$ and $(\beta ;\beta
+\varepsilon _0)$ do not contain any points from the support $S_{\,\varphi
}. $ Let $k_2$ be a natural number such that $r_{k_2}<\varepsilon _0$ and
let us consider set
\[
S_{\,k_2}=\left\{ x:\,x= \sum^{k_2}_{k=1} \gamma _k\,a_k,\,\gamma
_k\in \left\{ 0;1\right\} \right\} .
\]
It is easy to see that $S_{k_2}$ consists of at most $2^{k_2}$ elements. We
define $S_{k_2}^{^{\prime }}=S_{k_2}\cap [\alpha ;\beta ],$ and $%
S_{k_2}^{^{\prime \prime }}=S_{k_2}\backslash S_{k_2}^{^{\prime
}}.$ If a point $x_0=\sum\limits^{k_2}_{k=1} \gamma _k(x_0)\,a_k$
belongs to $S_{k_2}^{^{\prime \prime }},$ then any point $x$ of
the form
\[
x=\sum^{k_2}_{k=1} \gamma _k(x_0)a_k+ \sum^{\infty }_{k=k_2+1}
\gamma _k\,a_k
\]
does not belong to $[\alpha ;\beta ],$ since $x_0\leq \alpha -\varepsilon _0$
or $x_0\geq \beta +\varepsilon _0,$ and $r_{k_2}<\varepsilon _0.$

In the set $S_{k_2}^{^{\prime }}$ there exist a minimal point $x_1.$ It is
easy to understand that $x_1=\alpha $, otherwise the interval $[\alpha ;x_1)$
does not contain any points from $S_{\,\varphi }.$ Therefore
\begin{equation}
\alpha =\sum^{k_2}_{k=1} \gamma _k(\alpha )\,a_k+\sum^{\infty
}_{k=k_2+1} 0\cdot \,a_k=0,\gamma _1(\alpha )\gamma _2(\alpha
)\,\ldots \,\gamma _{k_2}(\alpha )\,000\ldots \label{alfa}
\end{equation}
Let $x_2$ be the second minimal point of the set $S_{k_2}^{^{\prime }},$
i.e.,
\[
x_2=\min \left\{ S_{k_2}^{^{\prime }}\backslash \,x_1\right\} .
\]
Finally, let us consider a natural number $k_3$ such that $%
k_3>k_2\,,\,a_{k_3}<(x_2-x_1)$ and $a_{k_3}>r_{k_3}.$ If $\gamma
_{k_2+1}(x)=1$ or $\gamma _{k_2+2}(x)=1$ or$\ldots $or $\gamma _{k_3}(x)=1,$
then $x\geq \alpha +a_{k_3};$ if $\gamma _{k_2+1}(x)=0$ and $\gamma
_{k_2+2}(x)=0$ and$\,\ldots $ and $\gamma _{k_3}(x)=0,$then $x\leq \alpha
+r_{k_3.}$ Therefore the interval $(\alpha +r_{k_3};\alpha +a_{k_3})$ does
not contain any point from the support. This \smallskip contradicts our
assumption.%
$\square$

\textbf{Remark.} The condition $a_k\geq a_{k+1},k\in N$ is not restrictive,
since the terms in the series (\ref{Series1}) are positive, and we can
choose the order of the terms of the series to fulfilled this condition.

Let us investigate criteria for singularity and absolutely continuity of the
distribution of the random variable $\varphi .$

\textbf{Theorem 12.} \textit{Let }$\delta _k>1$\textit{\ holds for} all $%
k\in N.\,$\textit{Then the distribution of the random variable }$\varphi $%
\textit{\ has pure type, and:}

\textit{1) is pure discrete if and only if }
\begin{equation}
\prod_{k=1}^{\infty } \max\limits_{0\leq i\leq s-1} \left\{
p_{0k},p_{1k}\right\}
>0; \label{discrAk}
\end{equation}
\textit{2) is pure absolutely continuous if and only if both}
\begin{equation}
\sum^{\infty }_{k=1} (\delta _k-1)<\infty ,\text{and }
\label{AbsDelta}
\end{equation}
\textit{\ }
\begin{equation}
\sum^{\infty }_{k=1}\,\,(\frac 12-p_{0k})^2\,<\infty \label{Sum1/2
Ak}
\end{equation}
converge;

\textit{3)\thinspace is pure singular in all other cases.\smallskip }

\textbf{Proof.}

\smallskip Purity of the distribution of the r.v. $\varphi $ follows from
the Jessen-Wintner theorem \cite{JW}.

Condition (\ref{AbsDelta}) is necessary for absolutely continuity because of
theorem 9. Let us consider two infinite product measure spaces $\left(
\Omega ,A,\mu \right) $and $\left( \Omega ,A,\nu \right) $ defined as in the
proof of theorem 7, i.e.,
\begin{eqnarray*}
\nu _k(0) &=&p_{0k},\,\,\,\nu _k(1)=p_{1k}; \\
\mu _k(0) &=&\frac 12,\,\,\,\,\,\,\mu _k(1)=\frac 12.
\end{eqnarray*}
It is easy to see that $\nu _k\ll \mu _k$ for all $k$, and $\rho (\mu _k,\nu
_k)=\int\limits_{\Omega _k} \sqrt{\frac{d\nu _k}{d\mu _k}}d\mu _k=\sqrt{%
\frac 12p_{0k}}+\sqrt{\frac 12p_{1k}}.$ By using theorem 2, we conclude that
the measure $\nu $ is absolute continuous with respect to measure $\mu \,$if
and only if
\begin{equation}
\prod^{\infty }_{k=1} \left( \sqrt{\frac 12p_{0k}}+%
\sqrt{\frac 12p_{1k}}\right) <\infty .  \label{InfProdAk}
\end{equation}
It is easy to check that the product (\ref{InfProdAk}) converges if and only
if condition (\ref{Sum1/2 Ak}) holds.

Let us consider a mapping $g$ : $\Omega \rightarrow [0;1]$ defined as
follows: for any $\omega =(\omega _1\omega _2\ldots \omega _k\ldots )\in
\Omega ,$ $\,$%
\begin{equation}
g(\omega )=\sum^{\infty }_{k=1} \omega _ka_k.
\end{equation}
It is easy to see that $g$ is a measurable mapping. Moreover, $g$ is
bijective, since $\delta _k>1$ holds for all $k\in N.$ The image measure $%
\mu ^{*}=\mu (g^{-1})\,$is a probability measure which is uniformly
distributed on the support $S_{\,\varphi }$ of positive Lebesgue measure $%
L(S_{\,\varphi })>0,$ i.e., for all Borel subset $E\subset [0;1]\,$we have
\[
\mu ^{*}(E)=\frac{L(E\cap S_{\,\varphi })}{L(S_{\,\varphi })}.
\]
The image measure $\nu ^{*}=\nu (g^{-1})$ is the distribution of the random
variable $\varphi .$ Since $g$ is a bijective mapping, the measure $\nu ^{*}$
is absolutely continuous with respect to the measure $\mu ^{*}$ if and only
if the measure $\nu $ is absolutely continuous with respect to the measure $%
\mu .$ Since $\mu ^{*}\ll L,$ we conclude that the measure $\nu ^{*}$ is
absolute continuous with respect to Lebesgue measure if and only if both (%
\ref{AbsDelta}) and (\ref{Sum1/2 Ak}) hold.

The criterium for discreteness follows from the Levi theorem \cite{Levi} or
can be deduced as a corollary from theorem 1 and the bijectivity of $g.$

Assertion 3) follows from the purity of the distribution of the r.v. $%
\varphi $. $\square$

\textbf{\smallskip Acknowledgment}

This work was partly supported by  DFG 436 UKR 113/53, INTAS 00-257,
SFB-611 projects and by DAAD. The last  named author gratefully acknowledged
the hospitality of the Institute of Applied Mathematics and of the IZKS of
the University of Bonn.

\end{document}